\documentclass{article}

    \usepackage[nonatbib,final]{neurips_2020_ml4ps}
   \usepackage{booktabs}

\usepackage[utf8]{inputenc} 
\usepackage[T1]{fontenc}    
\usepackage{hyperref}       
\usepackage{url}            
\usepackage{booktabs}       
\usepackage{amsfonts}       
\usepackage{nicefrac}       
\usepackage{microtype}      
\usepackage{graphicx}
\usepackage{amsmath}

\usepackage{caption}
\title{ADCME: Learning Spatially-varying Physical Fields using Deep Neural Networks}

\author{%
  Kailai Xu \\
  Institute for Computational and Mathematical Engineering\\
  Stanford University\\
  Stanford, CA 94305 \\
  \texttt{kailaix@stanford.edu} \\
   \And
  Eric Darve \\
  Institute for Computational and Mathematical Engineering and Mechanical Engineering\\
  Stanford University\\
  Stanford, CA 94305 \\
  \texttt{darve@stanford.edu} \\
}

\begin{document}

\maketitle

\begin{abstract}
  ADCME is a novel computational framework to solve inverse problems involving physical simulations and deep neural networks (DNNs). This paper benchmarks its capability to learn spatially-varying physical fields using DNNs. We demonstrate that our approach has superior accuracy compared to the discretization approach on a variety of problems, linear or nonlinear, static or dynamic. Technically, we formulate our inverse problem as a PDE-constrained optimization problem.  We express both the numerical simulations and DNNs using computational graphs and therefore, we can calculate the gradients using reverse-mode automatic differentiation. We apply a physics constrained learning algorithm (PCL) to efficiently back-propagate gradients through iterative solvers for nonlinear equations. The open source software which accompanies the present paper can be found at 
  \begin{center}
       \url{https://github.com/kailaix/ADCME.jl}
  \end{center}
\end{abstract}

\section{Introduction}

Inverse problems in computational engineering aim at learning physical parameters or spatially-varying fields from observations. These observations are usually the partial output of physical models, typically described by partial differential equations (PDEs). There is a vast array of literature covering inverse problems with an unknown spatially-varying field in geophysics \cite{landa2011joint,sen2013global,zhu2020general,li2020coupled}, fluid dynamics \cite{singh2016using,parish2016paradigm}, electromagnetism \cite{kelbert2014modem}, etc. In these applications, observations are indirect in the sense that they are not pointwise values of the unknown field. A standard approach is to discretize the spatially-varying field on the grid points, which is the same as our computational grid, and estimate pointwise values from observations. However, this approach is usually ill-posed when the dataset is small, which is common due to expensive and challenging experiments or measurements \cite{fomel2007shaping,huang2020learning}. 
 
Deep neural networks (DNN) are an effective approach to provide an expressive functional form for physical applications \cite{chandrasekhar2020tounn,hoyer2019neural}. It also provides some appropriate regularization to the optimization problem via proper initialization and choices of architectures \cite{raissi2019physics,hanin2018start,xu2019neural,xu2020learning}. In this paper, we present the capability of ADCME for learning spatially-varying physical parameters using deep neural networks \cite{ulyanov2018deep,rick2017one,fan2020solving}. In our approach, the inverse problem is formulated as a PDE-constrained optimization problem \cite{biegler2003large,rees2010optimal,de2015numerical,herzog2010algorithms}
\[
    \min_\theta \; J(u) \quad \text{($u$ is indirectly a function of $\theta$), \quad such that} \; F(u, \mathcal{N}_\theta) = 0
\]
Here $u$ is the state variable, $F$ is the PDE constraint, and $\mathcal{N}_\theta(\mathbf{x})$ is the approximation of the unknown field ($\mathbf{x}$ is the coordinate). 
The numerical solution $u$ can be formally expressed as a function of $\theta$ and we arrive at an unconstrained optimization problem by plugging $u(\theta)$ into $J(u)$ (see Figure 1)
\begin{equation}\label{equ:Ltheta}
    \min_\theta \; L(\theta) := J(u(\theta))
\end{equation}

\begin{center}
    \includegraphics[width=0.75\textwidth]{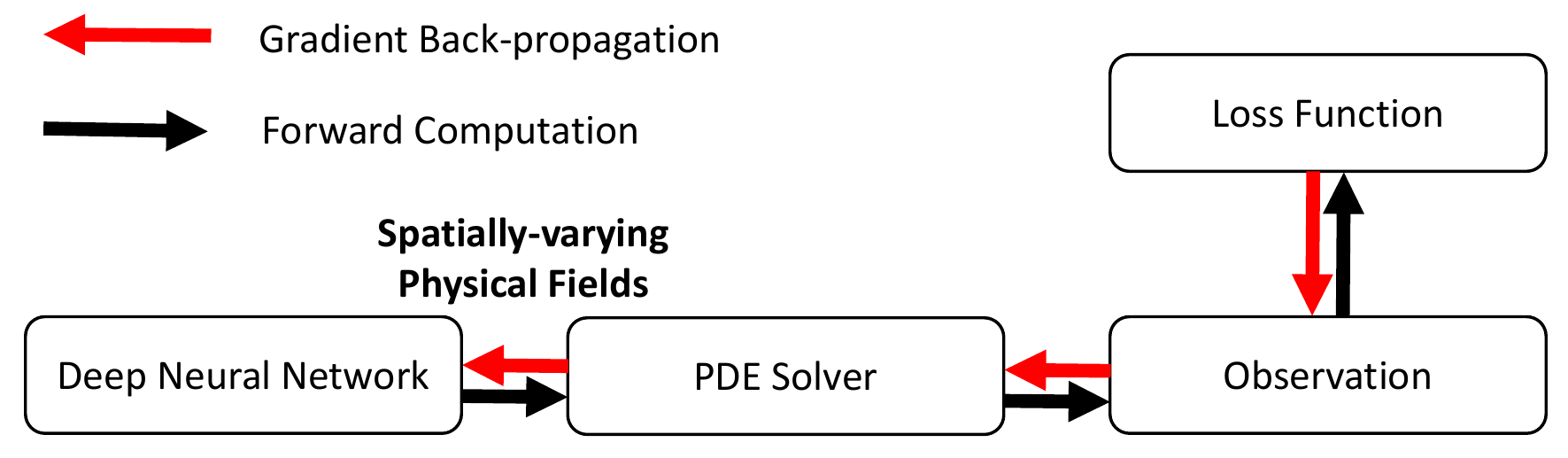}
    \captionof{figure}{Coupling a deep neural network and a numerical PDE solver.}
\end{center}
    

A gradient-based optimization algorithm is used to solve Eq.~\ref{equ:Ltheta}. The major challenge is to calculate $\nabla_\theta L(\theta)$. To this end, we express both numerical simulators (e.g., finite element methods) and deep neural networks as computational graphs. The advantage of this approach is that once we implement the forward computation, the gradient can be calculated using reverse-mode automatic differentiation (AD) \cite{baydin2017automatic,paszke2017automatic,abadi2016tensorflow}, which is equivalent to the so-called discrete adjoint methods \cite{plessix2006review,nachbagauer2015use,mcnamara2004fluid} for PDE solvers. ADCME provides a framework for implementing numerical schemes as computational graphs and incorporating deep neural networks. The computational graph operates on tensor operation (e.g., linear solvers) or higher levels (e.g., finite element assembling); with such task granularity, we enjoy a good balance between performance overhead from AD (e.g., reducing memory demands by checkpointing schemes for certain operators) and development efforts.   

Note that numerical simulators in ADCME are intrinsically different from traditional simulators: in ADCME, numerical simulators are defined by specifying the dependency of many individual and reusable operators, and these operators are equipped with gradient-backpropagation capabilities with the implementation of efficient adjoint state methods (possibly with checkpointing schemes for saving memories). 

With the advent of experimental techniques that enable gathering large amounts of data inexpensively, we anticipate that modeling using deep neural networks will become essential tools for data-driven discovery. As a result, developing software platforms with AD capabilities that can be coupled with sophisticated optimizers can be expected to have great impact on computational engineering.

\section{Mathematical Model}

Let us consider the Stokes problem as a concrete example. The governing equations for the Stokes problem are~\cite{mu2017simple}
\begin{gather}\label{equ:stokes}
    - \nabla \cdot (\nu \nabla u) + \nabla p = f \text{ in } \Omega \\ 
    \nabla \cdot u = 0 \text{ in } \Omega, \qquad
    u = 0 \text{ on }\partial \Omega \notag
\end{gather}
Here $u:\mathbb{R}^2\rightarrow \mathbb{R}^2$, $p:\mathbb{R}^2\rightarrow \mathbb{R}$ are state variables; $\nu(\mathbf{x})$ denotes the spatially-varying fluid viscosity, which we want to estimate from observations of $u$ or $p$; $\Omega\subset \mathbb{R}^2$ is the computational domain, and $f\in [L^2(\Omega)]^2$ is the unit external volumetric force. We approximate $\nu(\mathbf{x})$ with a deep neural network $\nu_\theta(\mathbf{x})$, which maps the coordinate $\mathbf{x}$ to a scalar value. $\theta$ are the DNN weights and biases. 

Another approach is to discretize $\nu(\mathbf{x})$ at the Gauss quadrature points or  in each finite element because only the values at Gauss quadrature points are used in the assembling procedures. Therefore, we have a vector of unknowns $\{\nu(\mathbf{x}_k)\}$, where $\{\mathbf{x}_k\}$ is the collection of all Gauss quadrature points of all the elements. To reduce the number of degrees of freedoms, we can also consider a piecewise constant function where $\nu$ is a constant on each element. 

The weak formulation of the Stokes equation finds $u \in [H^1_0(\Omega)]^2$ and $p\in L_0^2(\Omega)$ satisfying 
\begin{equation}
    \begin{aligned}
    (\nabla u, \nu_\theta\nabla v) - (\nabla \cdot v, p) &= (f, v), && v\in [H^1_0(\Omega)]^2 \\ 
    (\nabla \cdot u, q) &= 0, && q \in L_0^2(\Omega)
    \end{aligned}
\end{equation}
The finite element formulation using a P2/P0 (a quadratic velocity space and a pressure space of piecewise constant) element leads to the linear system 
$$\begin{bmatrix}A(\theta) & -B \\ B^T & 0\end{bmatrix} \begin{bmatrix}U \\ P\end{bmatrix}  = \begin{bmatrix}F_1 \\ F_2\end{bmatrix}$$
Note that the entries of $A$ are a function of $\theta$. In ADCME, constructing DNNs and assembling procedures of $A$, $B$, $F_1$, and $F_2$ are treated as nodes in the computational graph. These procedures are capable of back-propagating the gradients from downstream operators of the computational graph to upstream ones. When the forward computation is defined, a computational graph is constructed implicitly.  

The numerical scheme is used to compute the discrete solution $U$ and $P$, and we use $(U, P)$ to compute a loss function. Because we have constructed a computational graph, we are able to calculate the gradients by reverse-mode AD. One challenge here is that the numerical PDE solver may be nonlinear, leading to the solution of an implicit equation. To efficiently back-propagate the gradients, we apply the physics constrained learning method (PCL) proposed in \cite{xu2020physics} to extract the gradients. The core idea is to use the implicit function theorem to implement the adjoint rule~\cite{farrell2013automated,xu2020inverse}. 

\section{Numerical Experiments}

In this section, numerical examples, which cover linear and nonlinear elasticity, Stokes problems, and Burgers' equations, are presented to verify the effectiveness of our approach. We consider both a regular domain (square) and an irregular domain (a plate with two holes). The hyperelasticity problem involves nonlinearity in the numerical solver, and the Burgers' equation involves both nonlinearity and time stepping. As for numerical methods, the finite element method is used. For simplicity, the observations are all discrete solutions of the state variables unless specified. The neural network is a fully-connected deep neural network with 20 neurons per layer and 3 hidden layers. The activation function is \texttt{tanh}. We use the L-BFGS-B optimizer \cite{liu1989limited} with full batches. 

In what follows, we mention briefly our findings and the details of the numerical results can be found in the appendix. We will compare two approaches. The ``discretization'' methods assume that the values of $\nu(\mathbf{x}_k)$ can be optimized independently. Using our previous notations, we define $\theta = \{ \nu(\mathbf{x}_k) \}_k $ for all Gauss quadrature points $\mathbf{x}_k$. We also show benchmarks where $\nu(\mathbf{x})$ is approximated by a constant inside each element and we may write $\theta = \{ \nu(\mathbf{x}_e) \}_e$ where $e$ is now an element. The second approach represents the unknown function using a DNN.

We found that for discretization methods, the optimizer often leads to nonphysical solutions, i.e., values at some point are too large or negative, and this anomaly ends up breaking the numerical simulator. For DNNs, we do not impose extra constraints on the outputs except that we add a small positive values to the last layer to ensure DNNs start from a reasonable initial guess. Still, DNNs lead to a physical solution that is close to the exact solution. In summary, DNNs were found to be more stable and accurate than the discretization approaches in our problems.


\subsection{Linear Elasticity}

We consider the static linear elasticity problem on a unit square \cite{de2011computational}, where we have a spatially-varying Young's modulus $E(\mathbf{x})$ and a Poisson's ratio $\nu=0.3$. The governing equation for the linear elasticity problem is 
\begin{equation}\label{equ:elasticity}
\begin{aligned}
& \sigma_{ij,j}  +  b_i = 0, 
\quad
	\varepsilon_{ij} = \frac{1}{2}(u_{j,i}+u_{i,j}), \quad x\in \Omega\\
 & \sigma_{ij} = \lambda\delta_{ij} \varepsilon_{kk} + \mu(\varepsilon_{ij}+ \varepsilon_{ji}), \quad x\in \Omega\\
  &\sigma_{ij}n_j = t_j, \quad x\in \Gamma_N; \quad 	u_i = (u_0)_i, \ x\in \Gamma_D\\
&	\lambda = \frac{E\nu}{(1+\nu)(1-2\nu)},\quad
	\mu = \frac{E\nu}{1-\nu^2}
\end{aligned}
\end{equation}
Here $\sigma_{ij}$ is the stress tensor, $b_i$ is the body force, $u_i$ is the displacement, $\Gamma_N\cap \Gamma_D=\emptyset,\  \partial\Omega=\Gamma_N\cup \Gamma_D$. 

\begin{figure}[htpb]
    \centering
    \includegraphics[width=0.8\textwidth]{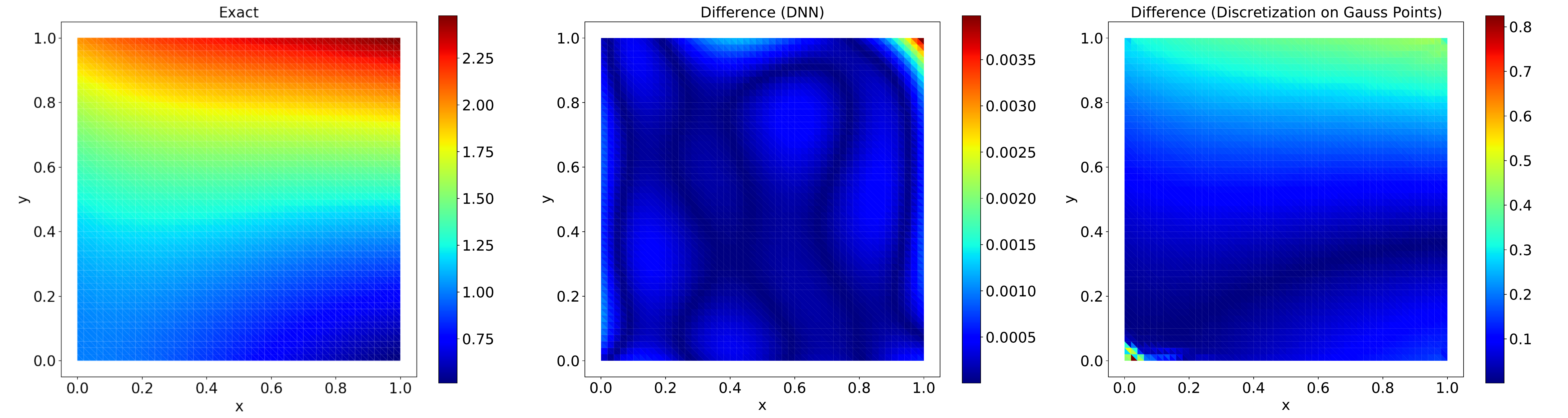}
    \caption{Convergence result for the static linear elasticity problem after 10,000 iterations.  For the discretization methods, we see we have some large errors at the lower left corner and top right region, and these anomalies break the numerical solver. The optimizer stops after 981 iterations. Note that \textbf{the scales for the second and third plots are different}.}
    \label{fig:linear_elasticity}
\end{figure}

We approximate $E(\mathbf{x})$ using a deep neural network. The results are shown in Figure~\ref{fig:linear_elasticity}. As a comparison, we also show results where we represent $E(\mathbf{x})$ as a discrete vector of trainable variables. We call this the ``discretization'' approach. We found that it did not work without imposing constraints on the trainable variables. In our experiment, the trainable variables were transformed from $E_i$ to $|E_i|$ (positivity constraint), where $E_i$ is the estimated $E(\mathbf{x})$ on the $i$-th Gauss points. Despite this, we can see that the DNN approach provides a better result that the discretization approach. 

\subsection{Stokes Problem}
The next problem considered is the Stokes problem Eq.~\ref{equ:stokes}.
The finite element formulation of the Stokes problem leads to a linear system, but has an extra pressure field compared to the linear elasticity problem. We assume that \textbf{only the pressure field is known} and we use this information to solve inverse problem. 

\begin{figure}[htpb]
    \centering
    \includegraphics[width=1.0\textwidth]{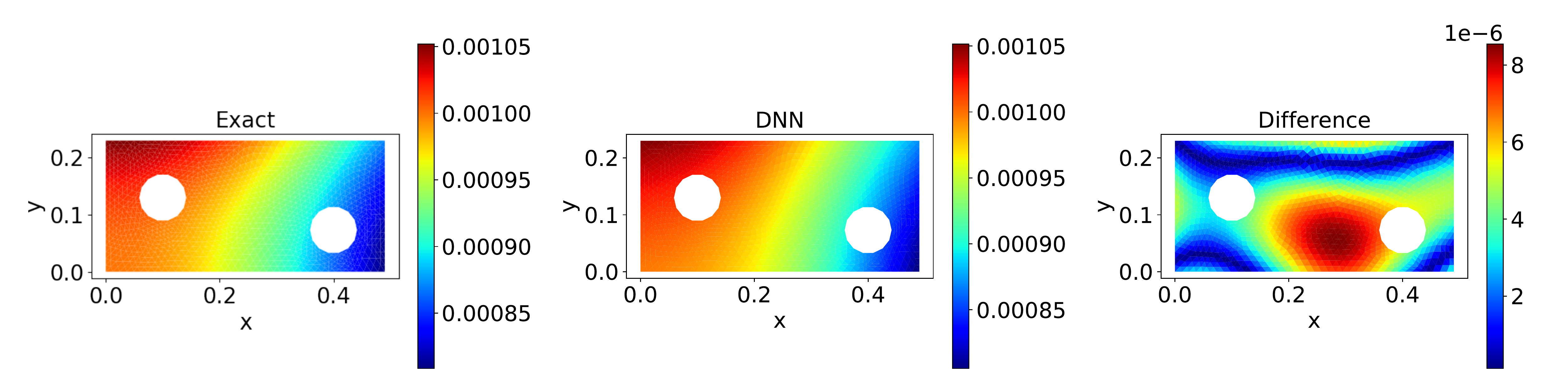}
    \caption{Convergence result for the Stokes problem after 600 iterations. For the discretization methods, the estimation did not converge to the correct results.}
    \label{fig:stokes}
\end{figure}

We make $\nu(\mathbf{x})$ a spatially-varying field. We were not able to obtain convergence result using the discretization approach even after imposing bound constraints on the trainable variables. This indicates that the inverse problem is quite ill-conditioned. The result for DNN is shown in Figure~\ref{fig:stokes}.

\subsection{Hyperelasticity}

The hyperelasticity problem can be described as an energy-minimization problem. We consider a common neo-Hookean  energy model~\cite{9Hyperel24:online} and assume that the Young's modulus $E(\mathbf{x})$ is a spatially-varying unknown field. 
The Lam\'e parameters are defined by the Young's modulus $E$ and Poisson ration $\nu$:
$$\lambda = \frac{E\nu}{(1+\nu)(1-2\nu)}, \quad \mu = \frac{E}{2(1+\nu)}$$
Using the finite element method to solve the energy-minimization problem leads to a nonlinear equation. The forward computation typically involves a Newton-Raphson iteration. To back-propagate the gradient through the nonlinear numerical solver, we applied the physics constrained learning to the Newton-Raphson solver. 

\begin{figure}[htpb]
    \centering
    \includegraphics[width=1.0\textwidth]{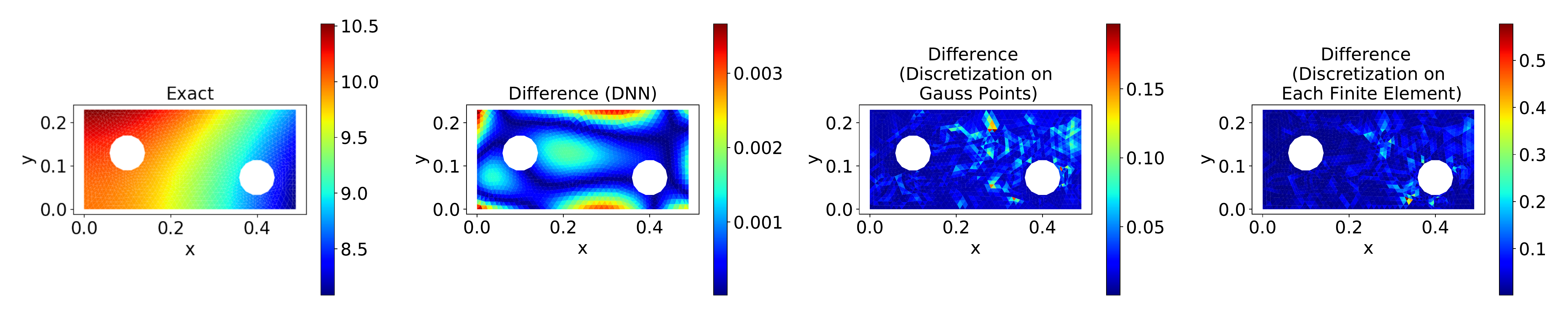}
    \caption{Convergence result for the hyperelasticity problem after 10,000 iterations. For the discretization method with piecewise constants, we use the prediction obtained at the 2,740th iteration, at which point the optimizer terminated.}
    \label{fig:hyperelasticity}
\end{figure}

We let the Young's modulus $E(\mathbf{x})$ be an unknown spatially-varying field. Figure~\ref{fig:hyperelasticity} shows the result. We see for this problem, the discretization approach provides an estimate that is close to the exact solution, and the corresponding loss function decays much faster than DNNs. However, the learned $E(\mathbf{x})$ field is not as smooth as DNNs and the maximum pointwise absolute error is much larger than DNNs. 


\subsection{Burgers' Equation}

The final example aims to demonstrate the capability of our algorithm for solving time-dependent nonlinear problems. We consider the Burgers' equation \cite{zhu2010numerical} with a spatially-varying viscosity parameter $\nu(\mathbf{x})$.

We use  the finite element method to discretize the Burgers' equation spatially and apply an implicit time integrator. This results in a time marching scheme, in which a nonlinear equation needs to be solved in each iteration. The physics constrained learning is used to solve the nonlinear equation. 

\begin{figure}[htpb]
    \centering
    \includegraphics[width=1.0\textwidth]{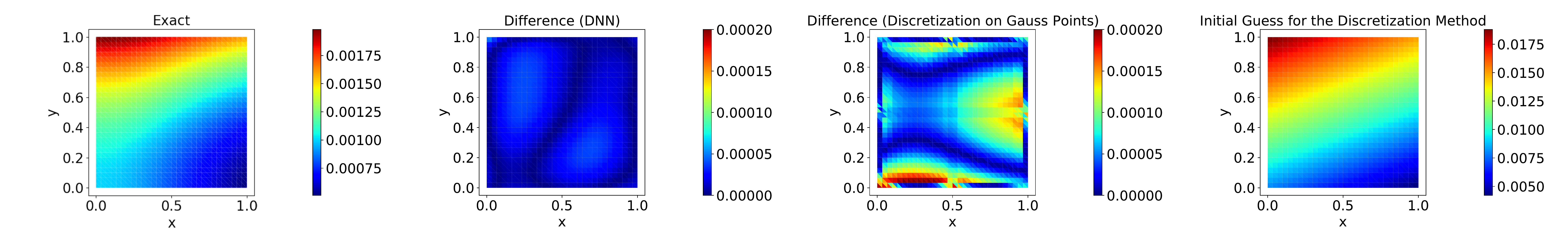}
    \caption{Convergence result for the Burgers' equation after 60 iterations (DNN) and 16 iterations (the Discretization Method). The optimization terminates at these iterations due to satisfying the stopping criterion.}
    \label{fig:burgers}
\end{figure}

We make $\nu(\mathbf{x})$ a spatially-varying field and apply both the DNN and the discretization approach to learn the physical field. The result is shown in Figure~\ref{fig:burgers}. Due to the ill-posedness of the inverse problem and lack of spatial correlation for the discretized $\nu$, we were not able to obtain convergence result with the discretization approach if we use a constant initial guess. To this end, we initialize the discretization using a field that is close to the exact field.\footnote{Other techniques, such as filtering \cite{sigmund200199}, can also be used for improving the discretization method.} Even with this impractical approach, the estimation is not as accurate as DNNs. 

\section{Conclusion}

We have presented the capability of ADCME for learning spatially-varying physical parameters using deep neural networks. We showed that the DNN approach gives us a much more accurate solution compared to the traditional discretization approach in the small data regime. In our approach, both numerical simulators and deep neural networks are expressed as computational graphs, and thus gradients can be back-propagated algorithmically once the forward computation is implemented. Particularly, for nonlinear solvers, we apply the physics constrained learning approach to back-propagate the gradients. The preliminary results on several PDE examples demonstrate the effectiveness of our approach. Note that the present work relies on an optimizer for a highly nonconvex problem and the well-posedness of the optimization problem, both of which remain elusive. We plan to investigate this direction further and develop a diagnosis toolbox in the future.

\newpage
\section*{Broader Impact}

This work allows developing more accurate mathematical models based on experimental data or based on a model developed at a finer scale (e.g., molecular scale). Advances in this paper will enable cheaper and more accurate computational models for a broad range of engineering and scientific applications. Computer simulations of engineering and physical systems are keys in many areas, including optimization, control, imaging, and failure analysis. Developing sophisticated computer models allows making accurate predictions for complex systems without having to resort to expensive experiments. Advances in sciences and engineering have a wide range of impacts on society, which are, for the most part, positive, such as transportation, energy, communication, and medicine, but, in some cases, may be harmful such as pollution and global warming or when applied to the development and production of weapons.

\section*{Acknowledgements}
This work is supported by the Applied Mathematics Program within the Department of Energy (DOE) Office of Advanced Scientific Computing Research (ASCR), through the Collaboratory on Mathematics and Physics-Informed Learning Machines for Multiscale and Multiphysics Problems Research Center (DE-SC0019453). 

\bibliographystyle{unsrt}
\bibliography{neurips}

\end{document}